\RequirePackage{fix-cm}
\documentclass[smallextended, numbook]{svjour3}       

\smartqed

\usepackage{a4wide,bm}
\usepackage{amssymb,verbatim}
\usepackage{euscript}
\usepackage{xcolor}
\usepackage{microtype}
\usepackage{graphicx}
\usepackage{amsmath}
\usepackage{subcaption}
\usepackage{float}
\usepackage{enumitem}   
\usepackage{mathrsfs} 
\usepackage[all]{xy}
\usepackage[ruled,vlined]{algorithm2e}
\usepackage[T1]{fontenc}
\usepackage{textcomp}
\usepackage{tikz}
\usepackage{capt-of}
\usepackage{mathtools}
\usepackage[colorlinks=true]{hyperref}
\hypersetup{
	colorlinks=true,
	linkcolor=red,
	filecolor=magenta,      
	citecolor=blue,
}
\usepackage[capitalize,noabbrev]{cleveref}

\newcommand{\nt}[1]{\textcolor{red}{#1}}

\DeclareMathOperator{\SL}{SL}
\DeclareMathOperator{\GL}{GL}
\DeclareMathOperator{\diag}{diag}
\DeclareMathOperator{\val}{val}

\newcommand{\R}{\mathbb{R}}
\newcommand{\K}{\mathbb{K}}
\newcommand{\F}{\mathbb{F}}
\newcommand{\Q}{\mathbb{Q}}
\newcommand{\Z}{\mathbb{Z}}

\newcommand{\Ocal}{\mathcal{O}}

\renewcommand{\P}{\mathbb{P}}

\newcommand{\Bcal}{\mathcal{B}}

\newcommand{\Tcal}{\mathcal{T}}

\newcommand{\Lcal}{\mathcal{L}}

\newcommand{\Xcal}{\mathcal{X}}
\newcommand{\Mcal}{\mathcal{M}}

\renewcommand{\char}{char}

\newcommand{\Prob}{\mathbb{P}}

\newcommand{\norm}[1]{\left\lVert#1\right\rVert}

\newcommand{\ind}{\text{ind}}

\newcommand{\indep}{\rotatebox[origin=c]{90}{$\models$}}

\begin{document}

\title{Statistics and tropicalization of local field Gaussian measures}

\author{Yassine El Maazouz         \and     Ngoc Mai Tran.
}

\institute{Yassine El Maazouz \at
          Department of Statistics, University of California at Berkeley, 335 Evans Hall \#3860, Berkeley,
          	CA 94720-3860, U.S.A.\\
              \email{yassine.el-maazouz@berkeley.edu}
           \and
           Ngoc Mai Tran \at
           Department of Mathematics, University of Texas at Austin, TX 78712 \\
           \email{ntran@math.utexas.edu}    
}

\date{Received: date / Accepted: date}

\maketitle

\begin{abstract} This paper aims to lay the foundations for statistics over local fields, such as the field of $p$-adic numbers. Over such fields, we give characterizations for maximum likelihood estimation and conditional independence for multivariate Gaussian distributions. We also give a bijection between the tropicalization of such Gaussian measures in dimension 2 and supermodular functions on the discrete cube $\{0,1\}^2$. Finally, we introduce the Bruhat-Tits building as a parameter space for Gaussian distributions and discuss their connections to conditional independence statements as an open problem.
 

\end{abstract}

\section{Introduction}
A local field $K$ is a locally compact, non-discrete  and totally discontinuous field. A typical example is the field of $p$-adic numbers $\Q_p$ for $p$ prime. There is an extensive literature on local fields in number theory \cite{cassels1986local}, analysis \cite{van1978non,schikhof2007ultrametric}, representation theory \cite{curtis1966representation}, mathematical physics \cite{vladimirov1994p,khrennikov2013p} and probability \cite{evans2001local,evans2007expectation,MR1822891}. While there has been an established theory for probability over local fields  \cite{evans2001local}, statistical problems such as maximum likelihood and conditional independence have not been considered to our best knowledge.  This paper aims to lay the foundations for this theory. 

One big motivator for us is the allure of a rich theory. There are nice analogues of fundamental objects in statistics over local fields, so one would hope that these objects are just nice enough to bring clean, complete characterizations that lead to an interesting and beautiful field.
Specifically, as shown by Evans \cite{evans2001local}, Gaussian distributions on $K^d$ display a tight link between orthogonality and independence. In general, $d$-dimensional $K$-Gaussians are parametrized by lattices in $K^d$. These lattices are analogous to the covariance matrix of real Gaussians, with the Bruhat-Tits building for the special linear group $\SL_d(K)$ plays the role of the cone of positive seimidefinite matrices. These properties allow solutions to statistical problems to be stated in terms of the geometry of the underlying lattices. In this paper we offer two such results, one on maximum likelihood estimation, the other on conditional independence of Gaussians on $K^d$. Definitions of relevant terms are given in the text.

\begin{theorem}\label{thm:mle}
	Let $ \Xcal = \{x_1, \dots , x_N\}$ be a dataset of points in $K^d$ of full rank. Then there is a unique full-dimensional Gaussian $G$ in $K^d$ that maximizes $\P(\Xcal|G)$, whose corresponding lattice is given by $L_{\Xcal} = \mathrm{span}(\Xcal)$.
\end{theorem}

\begin{theorem}\label{thm:ci}
	Let $X := (X_1,\dots, X_d)^T $ be a Gaussian vector in $K^d$ and $I$ a proper subset of $[d]$. The maximal subsets $J := \{ j_1, \dots , j_r \}$ of $\{1,\dots, d\} \setminus I$ such that $X_{j_1}, \dots X_{j_r}$ are mutually independent given $X_I$  are the bases of an $\F_q$-realizable matroid with base set $\{1,\dots, d\} \setminus I$ where $q$ is the module of $K$.
\end{theorem}

For simplicity we stated Theorem \ref{thm:mle} with the full-rank assumption. Generalization to the low rank case is discussed in Section \ref{sec:proof}. We also remark that the proof of Theorem~\ref{thm:ci} gives an explicit construction of the matroid. 

Our third theorem is motivated by the quest for the analogue of the Gaussian measure on the tropical affine space. In the recent years, this space has found fundamental applications in a diverse range of applications, from phylogenetics \cite{yoshida2017tropical,lin2018tropical} to social choice theory \cite{elsner2004max,tran2013pairwise}, game theory \cite{akian2012tropical} and economics \cite{baldwin2013tropical,tran2015product}. Of the various ways to define a `tropical Gaussian measure' \cite{tran2018tropical}, tropicalizing a Gaussian vector on a local field is the most theoretically attractive approach, for it opens up the possibility to formulate probabilistic questions in tropical algebraic geometry. We show that in dimension two, the tropicalization of such Gaussian measures is an interesting family of distributions which are in bijection with normalized supermodular function on the discrete cube $\{0,1\}^2$. We shall define the relevant terms in the text.

\begin{theorem}\label{thm:d2}
	Let $K$ be a local field with valuation $\val$ and module $q$. Let $X$ be a non-degenerate Gaussian vector in $K^2$ with lattice $L$ and $V:= \val(X)$ its image under valuation. Define $\varphi_L: \mathbb{Z}^2 \to \mathbb{R}$ via
	\[ \varphi_L(v) = - \log_q(\P(V \geq v)). \]
	Then, a function $\phi: \mathbb{Z}^2 \to \mathbb{R}$ equals $\varphi_L$ for some lattice $L$ if and only if $\phi$ is the restriction to $\Z^2$ of a tropical polynomial $P: \R^2 \to \R$ given by
	\begin{equation}\label{eqn:p.lambda}
	P(v) = \max( c_{00},  v_1 - c_{10},  v_2 - c_{01},  v_1 +  v_2 - c_{11}) ,
	\end{equation}
	where $c_{00} = 0$, and
	\begin{equation}\label{eqn:supermodular}
	c_{00} + c_{11} \geq c_{01} + c_{10}. 
	\end{equation}
	In other words, $P$ is supported on the discrete cube with supermodular coefficients.
\end{theorem}

Our paper is organized as follows. Section \ref{sec:background} reviews the essential background on Gaussian measures over local fields. Section \ref{subsec:mle} proves Theorem \ref{thm:mle} and \ref{thm:ci} and gives an algorithm to compute the defining lattice of a $K$-Gaussian. Section \ref{subsec:tropical} proves Theorem~\ref{thm:d2}. Section \ref{sec:buildings} discusses the structure of Bruhat-Tits buildings of the group $\SL_d(K)$. We conclude in Section \ref{sec:open} with discussions on two major research directions. The first concerns the relation between Bruhat-Tits buildings and conditional independence statements. The second is the generalization of Theorem \ref{thm:d2} to higher dimensions (cf. Example~\ref{ex:d3} and Conjecture~\ref{conj:d3}). We hope that this work will fuel more investigations in the novel area of statistics over local fields.

\section{Background and notations}\label{sec:background}

In this section we collect essential facts about Gaussian measures over local fields. Materials here are drawn from the monograph \cite{evans2001local} of Steve Evans, who established this area of study over a series of work \cite{evans1989local,evans1993local,evans2001local2,evans2002elementary,evans2007expectation,MR2266718}. For an extensive treatment on analysis over local fields, see \cite{van1978non,schikhof2007ultrametric}. We start with the fundamental example of Gaussians on the field of $p$-adic numbers in Section \ref{sec:background.p.adic}, and then give the general definitions and results in Section \ref{sec:background.general}.  

\subsection{Gaussians on  $\mathbb{Q}_p$}\label{sec:background.p.adic}

Fix a prime number $p$. We can write any non-zero rational number $r \in \Q \backslash \{0\}$ as $r = p^m a/b$ where $a,b \in \Z$ are integers not divisible by $p$ and $m \in \Z$ is unique. For example, with $p = 3$, then $\frac{8}{3} = 3^{-1} \times 8$, $5 = 3^0 \times 5$, $\frac{3}{8} = 3^{1} \times \frac{1}{8}$.  
We call $m$ the $p$-adic valuation of $r$ and write $\val_p(r) = m$.
The $p$-adic absolute value  $|\cdot|_p$ is
\[
|x|_p = p^{-\val_p(x)}, \quad x \in \Q.
\]
One can check that $|\cdot|_p$ is an \emph{ultrametric}. In particular, it satisfies the ultrametric inequality
\[
|x + y|_p \leq \max( |x|_p, |y|_p), \quad x,y \in \Q.
\]
The completion of $\Q$ with respect to the absolute value $|\cdot|_p$ is the field of $p$-adic numbers, denoted $\Q_p$. These numbers can be written as a Laurent series in $p$ with coefficients in $\Tcal = \{0, \dots, p-1\}$, that is
\begin{equation}\label{eqn:x.qp}
\Q_p =\left \{ x = \sum_{n \geq m } a_n p^n \colon m \in \Z, a_n \in \Tcal \right\}.    
\end{equation}
As noted in \cite[page 2]{evans2001local}, $\Q_p$ has a fractal-like structure. It can be visualized as a rooted $p$-ary tree. The root of this tree is the valuation ring $\Z_p$ (the $p$-adic integers)
\begin{equation}\label{eqn:x.zp}
\Z_p = \{x \in \Q_p, x = \sum_{n \geq 0 } a_n p^n, \quad  a_n \in \Tcal\}. 
\end{equation}
Note that $\Z_p = \bigcup_{a \in \Tcal}(a + p\Z_p)$, that is, it is the union of $p$ disjoint translates of itself. These $p$ translates are the $p$ children of $Z_p$ at level 1, and by recursion, each child in turn has $p$ children made up of $p$ translated of copies of itself, and so on. We can extend backwards to obtain a tree with levels indexed by $\mathbb{Z}$. For a fixed $x \in \Q_p$, the coefficients $(a_n \in \Tcal)$ in its expansion \eqref{eqn:x.qp} tells us the sequence of nodes of the tree that $x$ belongs to. 

\begin{example} $\Q_p$  is an infinite tree with degree $p+1$. \cref{fig:Q_2} is a local depiction of $\Q_p$.
	
	\begin{figure}[!ht]
		\centering
		\resizebox{1\textwidth}{!}{
			\tikzset{every picture/.style={line width=0.75pt}}      
			
			\begin{tikzpicture}[x=0.5pt,y=0.6pt,yscale=-1,xscale=1]
			\draw    (321,83) -- (130,129) ;
			\draw    (501,132) -- (321,83) ;
			\draw    (116,156) -- (36,237) ;
			\draw    (198,238) -- (116,156) ;
			\draw    (500,153) -- (421,230) ;
			\draw    (582,229) -- (500,153) ;
			
			\draw (313,60.4) node [anchor=north west][inner sep=2pt]    {$\mathbb{Z}_{2}$};
			\draw (97,136.4) node [anchor=north west][inner sep=2pt]    {$2\ \mathbb{Z}_{2}$};
			\draw (473,135.4) node [anchor=north west][inner sep=2pt]    {$1\ +\ 2\mathbb{Z}_{2}$};
			\draw (15,242.4) node [anchor=north west][inner sep=2pt]    {$2^{2} \ \mathbb{Z}_{2}$};
			\draw (161,242.4) node [anchor=north west][inner sep=2pt]    {$2+2^{2} \ \mathbb{Z}_{2}$};
			\draw (388,232.4) node [anchor=north west][inner sep=2pt]    {$1+2^{2} \ \mathbb{Z}_{2}$};
			\draw (545,231.4) node [anchor=north west][inner sep=2pt]    {$1+2+2^{2} \ \mathbb{Z}_{2}$};
			
			\end{tikzpicture}
			
		}
		\caption{Local depiction of $\Q_2$ as a tree.}
		\label{fig:Q_2}
	\end{figure}
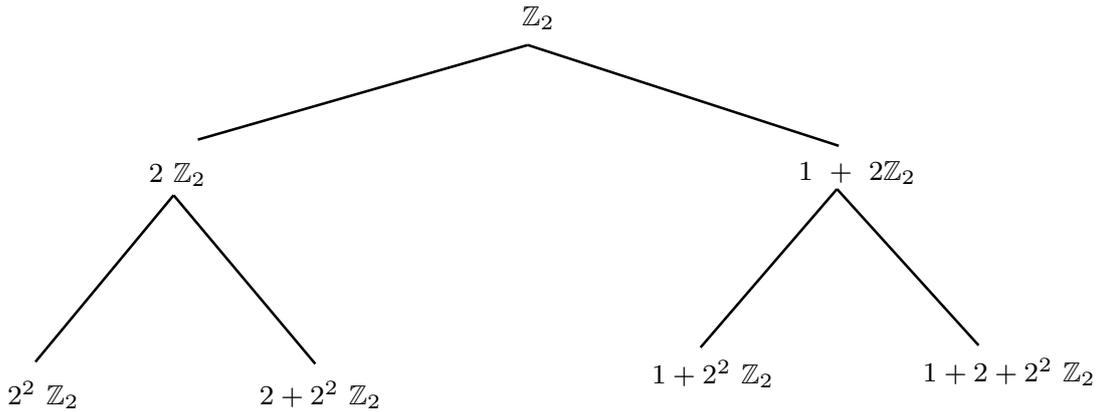
	
\end{example}

Parallel to Kac's characterization of classical Gaussians \cite{kac1939characterization}, Evans \cite[Definition 4.1]{evans2001local} defined the Gaussian measure on $K^d$ for some local field $K$ to be one that is invariant under orthonormal transformations. 

\begin{definition}[\cite{evans2001local}, Definition 4.1]\label{def:Gaussians}
	A random variable $X$ on $K^d$ has a centered Gaussian distribution if whenever $X_1,X_2$ are independent copies of $X$, and $A = \begin{bmatrix} a_{11} & a_{12} \\ a_{21} & a_{22} \end{bmatrix}$ is a matrix in $GL_2(\Ocal_K)$ then 
	\[ 
	\begin{bmatrix} X_1 \\ X_2 \end{bmatrix} \stackrel{d}{=} A \begin{bmatrix} X_1 \\ X_2 \end{bmatrix} := \begin{bmatrix} a_{11} X_1 + a_{12} X_2 \\ a_{21} X_1 + a_{22} X_2 \end{bmatrix} 
	\]
\end{definition}

For $K = \Q_p$ and $d = 1$, \cite[Theorem 4.2]{evans2001local} says that each Gaussian in $\Q_d$ is the uniform probability measure supported on $p^n \Z_p$ for some $n \in \Z$. For $d > 1$, the story is much more interesting. Namely, the Gaussian distributions on $\Q_p^d$ are exactly the uniform probability measures on $\Z_p$-submodules of $\Q_p^d$. These are set of the form
\[
L = \Z_p u_1 + \dots + \Z_p u_n = \left\{ \sum_{i = 1}^{d} a_i u_i \colon a_i \in \Z_p \text{ and } u_i \in \Q_p^d  \right\}.
\]
Visually, $L$ is a lattice in $\Q_p^d$ with generators $\{u_1,\dots,u_d\}$. Each lattice corresponds to the support of a unique Gaussian measure, and vice versa. 
\begin{example}
	The uniform probability measure on $\Z_3^2$ is the \emph{standard Gaussian} measure on $\Q_3^2$. Its lattice $L$ is the standard lattice generated by the standard vectors $u_1 = (1,0)$, $u_2 = (0,1)$. 
	If $X = (X_1, X_2)$ has this probability distribution, then coordinates $X_1,X_2$ are independent standard Gaussians on $\Z_3$. The vector $Y = AX$ where 
	\[
	A = \begin{bmatrix}
	1 & 1 \\ 3 & 0
	\end{bmatrix},
	\]
	is a Gaussian vector with uniform distribution on the module
	\[
	L =  \Z_3 \begin{bmatrix} 1 \\ 3 \end{bmatrix}  + \Z_3 \begin{bmatrix} 1 \\ 0 \end{bmatrix}.
	\]
	The coordinate $Y_1, Y_2$ are Gaussians on $\Q_3$ but they are no longer independent.
\end{example}

The lattice $L$ here plays the role of the covariance matrix of classical Gaussians over $\mathbb{R}$. The reader might be wary with the choice of the matrix $A$ not being unique, but the situation however is very much similar to the real case since, given a non-degenerate covariance matrix $\Sigma \in \R^{d \times d}$ the matrix equation $B B^\top = \Sigma$ has an infinite number of solutions, namely the left coset $B O_d(\R)$ where $O_d(\R)$ is the orthogonal group and $B$ is an arbitrary solution. The simlimarity goes even further because, given a lattice $L$ in $\Q_p^d$, there exists a matrix $A \in \Q_p^{d \times d}$ such that $L = A \Z_p^d$ and the solutions for this equation are exactly the matrices in the coset $A \GL_d(\Z_p)$. The group $\GL_d(\Z_p)$ plays the role of the orthogonal group in this setting as we shall see below (cf. Corollary \ref{cor:gl.d}). 

\subsection{General local field}\label{sec:background.general} 
This section generalizes the above discussion from $\Q_p$ to a non-archimedean local field $K$, that is, a locally compact, non-discrete, totally disconnected, topological field. We denote by $K^\times$ be the set of all invertible elements in $K$. The field $K$ comes with an additive, surjective, discrete valuation map $\val : K \rightarrow \Z  \cup \{+ \infty \}$. We fix an element $\pi \in K$ such that $\val(\pi) = 1$, such an element is called a \emph{uniformizer of} $K$. The \emph{valuation ring} of $K$ is $\Ocal_K := \{ x \in K \colon  \val(x) \geq 0 \}$. This is a discrete valuation ring with the unique maximal ideal $\pi \Ocal_K = \{x \in \Ocal_K \colon: \val(x) > 0 \}$. The \emph{residue field} $k = \Ocal_K / \pi \Ocal_K$ is isomorphic to a finite field $\F_q$ of cardinality $q = p^c$ where $p :=\char(k)$ is prime and $c \geq 1$ is an integer. The cardinality $q$ of the residue field $k$ is also called the \emph{module} of $K$. The field $K$ is equipped with an absolute value $|\cdot|$ defined as
\[
|x| = q ^{ - \val(x)}, \quad x \in K.
\]

\begin{example}
	For $K = \Q_p$, $\val$ is the $p$-adic valuation, $\pi =p$, $\Ocal_K = \Z_p$ and $q = p$. 
\end{example}

The following is the generalization of \eqref{eqn:x.qp}. 
\begin{proposition}[\cite{ash2003numbertheory}, \S 9.4.4 ]\label{prop:series}
	Fix a set $\Tcal \subset \Ocal_K$ of representatives of elements in $k = \Ocal_K / \pi \Ocal_K$ such that $\Tcal$ contains $0$. Let $x \in K^{\times}$. There exists a unique integer $n \in \Z$ and a unique sequence $(u_i)_{i \geq n} $ of elements in $\Tcal$ such that:
	
	$$ x = \sum \limits_{i \geq n} u_i \pi^i $$
\end{proposition}

Now we collect results concerning lattices in $K^d$. 
Let $d \geq 1$ be a positive integer. For the \emph{standard lattice} $\Ocal_K^d$, define the norm $\norm{\cdot}$ by
\[
\norm{x} := \max\{ |x_i| \colon x = (x_1, \dots, x_d) \in K^d \}.
\]
For $x \in \Ocal_K$, define $\overline{x} := (x \mod \pi) \in k$. Similarly, if $x = (x_1, \dots, x_d)$ is a vector in $\Ocal_K^d$, we denote by $\overline{x} = (\overline{x_1}, \dots, \overline{x_d}) \in k^d$. 

We denote by $K^{\vee d}$ the dual space of $K^d$ (the space of linear forms on $K^d$). There exists a natural definition of orthogonality on $K^d$ given as an analogue for Pythagoras' theorem in Euclidean spaces. As in the classical settings, orthogonality implies linear independence. 

\begin{definition}\label{defn:orthogonality}
	Let $F = \{x_1, \dots, x_n\} \subset K^d$ be a collection of vectors. We say that $F$ is orthogonal if 
	\[ 
	\norm{\sum \limits _{i=1}^{n} \alpha_i x_i} = \max\limits_{1\leq i \leq n} |\alpha_i|\norm{x_i}, \quad \text{ for any }  \alpha_1, \dots, \alpha_n \in K.
	\]
	If, in addition, the vectors $x_1, \dots, x_n$ all have norm $1$ the set $F$ is called orthonormal.
\end{definition}

Note that the orthogonality of a set of non-zero vectors $ x_1, \dots, x_n $ in $K^d$ does not change if we scale the vectors with non-zero scalers in $K^\times$. The following proposition gives a more practical criterion of orthogonality.

\begin{proposition}[\cite{van1978non}, Exercise 5.A] \label{prop:orthogonality} 
	Let $x_1, \dots, x_n \subset K^d $ be a finite set of vectors in $K^d$ of norm $1$. In particular, $x_i \in \Ocal_K^d$ for all $i = 1, \dots, d$. Then
	$x_1, \dots, x_n$ are orthonormal if and only if  $\overline{x_1}, \dots, \overline{x_n}$ are linearly independent in $k^d$.
\end{proposition}

A $d \times d$ matrix in $\Ocal_K^{d \times d}$ is said to be orthogonal if its row vectors are orthonormal. The set of such orthogonal matrices has a group structure as stated by the following corollary. An important consequence is that there is an analogue of a Gram-Schmidt process on $K^d$, and a singular value decomposition for matrices over $K$. 

\begin{corollary}\label{cor:gl.d}
	The set of $d \times d$ matrices with orthonormal columns in $K^d$ is exactly the group $\GL_d(\Ocal_K) =\{ U \in \GL_d(K)\cap \Ocal_K^{d \times d} : U^{-1} \in \Ocal_K^{d \times d} \} $ of invertible matrices in $\Ocal_K^{d \times d}$.
\end{corollary}

\begin{proposition}[\cite{schikhof1984ultrametric}, Theorem 50.8]\label{prop:gram.schmidt}
	Let $(e_k)_{1\leq k \leq n}$ be a collection of linearly independent vectors in $K^d$. There exists a collection of orthonormal vectors $(v_k)_{1\leq k \leq n}$ such that $\mathrm{span}_K(  (e_i)_{1 \leq i \leq k} ) =  \mathrm{span}_K(  (v_i)_{1 \leq i \leq k} ) $ for all $k \in \{1,\dots, n\} $.
\end{proposition}

The following proposition is the nonarchimedean analog of singular value decomposition (SVD). Since the group of orthogonal matrices is exactly $\GL_d(\Ocal_K)$, the Smith normal form and SVD are the same concept.

\begin{proposition}[\cite{evans2002elementary}, Theorem 3.1]\label{prop:non_arch-SVD}
	Let $A \in K^{d \times n}$, there exist two orthogonal matrices $U \in \GL_d(R)$ and $V\in \GL_n(R)$ and a matrix $D \in K^{d \times n}$ such that $A = UDV$ and all off diagonal entries of $D$ are zero.
\end{proposition}

A \emph{lattice} in $K^d$ is a compact $\Ocal_K$-submodule of $K^d$. By \cite[Chapter II-Proposition 5]{weil2013basic}, all compact $\Ocal_K$-submodules of $K^d$ are finitely generated. Thus, lattices of $K^d$ are of the form $A \Ocal_K^m$ for some $A \in K^{d \times m}$. In this paper, we shall frequently use one of two canonical choices for $A$:  the orthonormal form and the Hermite normal form (cf. Definition \ref{defn:lower.triangular}). The existence and uniqueness of these forms follow from row-reduction operations over $K$, analogous to classical proofs over $\mathbb{R}$.

\begin{lemma}[Orthonormal form of a lattice] \label{Lem:normal_form}
	Let $n \leq d$ be an integer. Lattices of rank $n$ in $K^d$ are exactly those of the form $\pi^{m_1} \Ocal_K u_1 \oplus \dots \oplus \pi^{m_n} \Ocal_K u_n$ for some sequence of integers $m_1, \dots , m_n$ and orthonormal vectors $u_1, \dots, u_n$ of $K^d$.
\end{lemma}

\begin{proof}
	Follows immediately from Proposition \ref{prop:non_arch-SVD}.
\end{proof}

\begin{definition}\label{defn:lower.triangular}
	Say that a matrix $A = (a_{ij})_{1 \leq i,j \leq d} \in \GL_d(K)$ in Hermite normal form if 
	
	\begin{enumerate}
		\item $a_{ij} = 0$  for all $1 \leq i < j \leq d$,
		\item $a_{ii} = \pi^{n_i}$ with $n_i \in \Z$ for all $1 \leq i \leq d$,
		\item $a_{ij}$ is either $0$ or is Laurent polynomial in $\pi$ with coefficients in $\Tcal$ of degree less than $n_i$ for all indices $1 \leq j < i \leq d$.
	\end{enumerate}
	
\end{definition}

\begin{example}\label{example:HermiteFormMatrices}
	Consider the two following matrices in$\Q_2^{3 \times 3}$
	\[
	\begin{bmatrix}
	1 & 0 & 0\\ 1/2 & 2 &0 \\ 1& 2 & 2^2   
	\end{bmatrix}  
	\quad \text{ and } \quad
	\begin{bmatrix}
	1 & 0 & 0\\ 1/2 + 2^2 & 2 &0 \\ 2^2  & 2 & 1
	\end{bmatrix}.
	\]
	The first matrix in $\Q_2^{3 \times 3}$ is in Hermite normal form but the second is not.
	
\end{example}

\begin{lemma}[Hermite normal form a lattice]\label{lem:Normal_Form}
	Let  $L$ be a lattice of rank $d$. There exists a unique matrix $A = (a_{ij}) \in \GL_d(K)$ such that $A$ is in Hermite form, and $ L = A \Ocal_K^d$. We call $A$ the Hermite normal form of the lattice $L$. 
\end{lemma}

\begin{proof}
	We start by proving uniqueness. Let $A$ and $B$ be two matrices satisfying the conditions above. Since $A$ and $B$ represent the same lattice, each column vector of $A$ is a linear combination with coefficients in $\Ocal_K$ of the columns of $B$ and vice versa. Thanks to the lower triangular form of $A$ and $B$, the $i$-th column of $A$ is a linear combination of the columns indexed by $i, i+1, \dots , d$ in $B$ with coefficients in $\Ocal_K$. By Definition \ref{defn:lower.triangular}, it follows  that $A$ and $B$ have the same last column. Let $ 1\leq j \leq d-1$ and suppose that columns indexed by $j+1, \dots, d$ in $A$ and $B$ are identical. We have again $A_{jj} = B_{jj}$ thanks to the lower triangular form of $A$ and $B$. The condition on the power series representation of $A_{ij}$ for $i>j$ of Definition \ref{defn:lower.triangular} allows us to conclude that the $j^{th}$ column of $A$ and $B$ are equal. Thus $A = B$ which settles uniqueness.
	To prove existence, it suffices to  transform $A$ to a  matrix $A'$ in Hermite form by multiplying on the right with elementary matrices (with entries in $\Ocal_K$) and permutation matrices which are all elements of $\GL_d(\Ocal_K)$ i.e. we want to find a matrix $U \in \GL_d(\Ocal_K)$ such that $AU$ is in Hermite normal form. Since orthogonal matrices stabilize the standard lattice $\Ocal_K^d$ we shall get $L = A U \Ocal_K^d = A' \Ocal_K^d$. We now explain how to get the matrix $U$. We may assume that $\val(a_{11}) \leq \val(a_{1j})$ for all $1 \leq j \leq d$, otherwise we can multiply the matrix $A$ with a permutation matrix on the right to permute the columns. Multiplying $A$ on the right by the matrices elementary $I - \frac{a_{1j}}{a_{11}} E_{1j} \in \GL_d(\Ocal_K)$ for $1 < j \leq d$, cancels the entries $(1,j)$ for $1 < j \leq  d$. Repeating this process with the remaining $d-1$ columns we can choose $A = (a_{ij})$ to be lower triangular. Multiplying $A$ on the right with the diagonal matrix $\diag(\frac{\pi^{\val(a_{11})}}{a_{11}}, \dots, \frac{\pi^{\val(a_{dd})}}{a_{dd}}) \in GL_d(\Ocal_K)$ we can choose $A$ such that conditions (1) and (2) are satisfied. All that remains is to satisfy condition (3). Using \cref{prop:series} and multiplying with elementary matrices of the form $I -  \alpha E_{ij} \in \GL_d(\Ocal_K)$ with $1 \leq j < i \leq d$ and $\alpha \in \Ocal_K$ we can satisfy condition (3).
\end{proof}

\begin{example}
	Let $A$ be the second matrix in \cref{example:HermiteFormMatrices} and $L = A \Z_2^3 \subset \Q_2^3$ its corresponding lattice. If we multiply $A$ on the right by the matrix
	\[
	U = \begin{bmatrix}
	1 & 0 &0 \\	  -2 &1 &0\\ -2^2 & -2 &1 
	\end{bmatrix} \in \GL_3(\Z_2),
	\]
	we get the matrix 
	\[
	AU =	\begin{bmatrix}
	1 & 0 & 0\\ 1/2 & 2 &0 \\ 0  & 0 & 1
	\end{bmatrix},
	\]
	which is in Hermite normal form and represents the same lattice $L$.
\end{example}

Finally, we now recall important results that connect independence of Gaussians with orthogonality. We denote by $\lambda$ the unique Haar measure on $K$ such that $\lambda(\Ocal_K) = 1$ \cite[p4]{evans2001local}. The space $K^d$ is then equipped with the product measure induced by $\lambda$. With no risk of confusion, we also denote the product measure on $K^d$ by $\lambda$. 

Recall the definition of Gaussians on $K^d$ (cf. \cref{def:Gaussians}). Evans completely characterized all non-trivial Gaussian measures on $K^d$. The following is a rephrasing of his results (Theorems 4.4 and 4.6 of \cite{evans2001local}) in terms of lattices over~$K^d$. 

\begin{theorem}\label{prop:multivar_gauss}
	The distributions of $K^d$-valued Gaussian random variables are exactly the normalized Haar measures on lattices in $K^d$.
\end{theorem}

In other words, the law of each Gaussian on $K^d$ is completely specified by a lattice $L$. This is analogous to the classical case, where the law of each centered Gaussian in $\mathbb{R}^d$ is completely specified by its covariance. We note that  for Gaussians over a local field the mean is not well-defined \cite{evans2001local}, so lattices in $K^d$ are indeed the central objects of the theory of Gaussian over $K$. We say that a Gaussian measure on $K^d$ is \emph{non-degenerate} if its corresponding lattice has full rank. We call \emph{standard Gaussian distribution} on $K^d$ the Gaussian distribution defined by the standard lattice $\Ocal_K^d$ i.e. the uniform probability measure on $\Ocal_K^d$ with respect the Haar measure on $K^d$. As in the Euclidean setting, independence of Gaussians is tightly linked to orthogonality as the following Lemma explains.

\begin{lemma}[\cite{evans2001local}, Theorem 4.8] \label{lem:indep}
	Let $f_1,...,f_n \in K^{\vee d}$ be linear forms (identified as row vectors of a matrix) and $Z$ a standard Gaussian vector in $K^d$. Then,  $(f_i(Z))_{1 \leq i \leq n}$   are mutually independent if and only if  $f_1, \dots, f_n$ are orthogonal.
\end{lemma}

\begin{example}
	Suppose $d=4$ and let $A \in \Q_7^{4 \times 4}$ be the following matrix 
	\[
	A = \begin{bmatrix} 12 \ & 314 \ &  234 \ & 34  \\ 12 \ & 343  \ & 55 \  & 67   \\ 25 \  & 54 \  &  65 \  &  65  \\61 \  & 461 \  & 430 \  & 328  \\   \end{bmatrix}.
	\]
	Let us use \cref{prop:orthogonality} to test the orthogonality of the rows of $A$. Every row in this matrix has norm $1$ in $\Q_7^4$ and modulo $7$ the matrix becomes
	\[
	\overline{A} = \begin{bmatrix} 	5 \ & 6 \ & 3 \ & 6  \\ 5 \  & 0 \ & 6 \  & 4 \ \\ 4 \ & 5 \ &  2 \ & 2  \\ 5 \ & 6 \ & 3 \ & 6  \end{bmatrix} \in \F_7^{4 \times 4}.
	\]
	So if $Z = (Z_1, Z_2, Z_3,Z_4)^{\top}$ is a vector of independent standard Gaussians in $\Q_7$ and $Y = AZ$ we can see by \cref{lem:indep} that $Y_1,Y_2,Y_3$ are independent but $Y_1,Y_2,Y_3,Y_4$ are not since $\overline{A}\in \F_7^{4 \times 4}$ only has rank $3$.
\end{example}

\section{Proof of main results}\label{sec:proof}
\subsection{Maximum likelihood and conditional independence}\label{subsec:mle}
We now prove Theorems~\ref{thm:mle} and~\ref{thm:ci} stated in the introduction. Given a full rank lattice $L \subset K^d$, normalized Haar measure on $L$ is a Gaussian distribution and it has a probability density function $ f_L = \frac{1}{\lambda(L)} 1_L$. Let $ \Xcal = \{ x_1, \dots , x_N \}$ be a dataset of points in $K^d$ of full rank. The likelihood of observing $\Xcal$ assuming that the data coming from the Gaussian distribution with lattice $L$ is given by
$$\mathscr{L} (\Xcal,L) = \prod \limits _ {i=1}^{N} f_L(x_i).$$ 
Our goal in the settings of Theorem \ref{thm:mle} is to maximize the likelihood in terms of $L$.
Note that this is equivalent to finding  a lattice $L$ that contains all the data points $\Xcal$ and has minimal $\lambda$ measure. The proof relies on the following result, which gives an expression for $\lambda(L)$ in terms of its matrix representation. 

\begin{lemma}\label{lem:measure_exp}
	Let  $L = A \Ocal_K^d$ be a rank $d$ lattice in $K^d$ where $A \in \GL_d(K)$. Then we have 
	\[ 
	\lambda(L) = q^{-\val(\det(A))} . 
	\]
\end{lemma}
\begin{proof}
	By Proposition \ref{prop:non_arch-SVD}, there exists an orthogonal matrix $U \in \GL_d(\Ocal_K)$ and a sequence $n_1, \dots, n_d$ of integers such that $L = U \diag(\pi^{n_1}, \dots, \pi^{n_d}) \Ocal_K^d$. We then have 
	\[
	\diag(\pi^{n_1}, \dots, \pi^{n_d}) \Ocal_K^d = \pi^{n_1} \Ocal_K \times \dots \times \pi^{n_d} \Ocal_K
	\]
	and then $\lambda(\diag(\pi^{n_1}, \dots, \pi^{n_d}) \Ocal_K^d ) = \prod \limits _{i = 1} ^{d}\lambda(\pi^{n_i}\Ocal_K) = q^{- (n_1 + \dots + n_d)}$. Since orthogonal matrices preserve the measure we also have $\lambda(L) = q^{- (n_1 + \dots + n_d)} $. Furthermore, since $\val(\det(V)) = 0$ for all $V \in \GL_d(\Ocal_K)$ we have $\val(\det(A)) = n_1 + \dots + n_d$.
\end{proof}

\begin{proof}[Proof of Theorem \ref{thm:mle}]
	Define the lattice $L_{\Xcal} := \mathrm{span} _{\Ocal_K}(\Xcal)$. Since $\Xcal$ is of full rank, $L_{\Xcal}$ is also full rank. Now let $L$ be any other lattice that contains $\Xcal$. Then $\Lcal(\Xcal,L) = \lambda(L) ^ {-1}$ and $ \Lcal(\Xcal,L_{\Xcal}) = \lambda(L_{\Xcal}) ^ {-1}$ . Since $L_{\Xcal} \subset L$ we have $\lambda(L) \geq \lambda(L_{\Xcal})$. Thus $\Lcal(\Xcal,L) \leq \Lcal(\Xcal,L_{\Xcal})$. The lattice $L_{\Xcal}$ maximizes the likelihood. Suppose that $L_{\Xcal} \subsetneq L$ by means of a basis change without loss of generality we can suppose that $L = \Ocal_K^d$. There exists an orthogonal matrix $U$ and a diagonal matrix $D := diag(\pi^{n_1},\dots , \pi^{n_d})$ such that $L_{\Xcal} = U D \Ocal_K^d $. Since $L_{\Xcal} \subsetneq L$ we have $n_i \geq 0$ for any $1 \leq i \leq d $ and there exists $k$ such that $n_k > 0$. Thus $\lambda(L_X) < \lambda(L)$. Then $L_X$ is the unique lattice that maximizes the likelihood. \end{proof}

\begin{proof}[Proof of Theorem \ref{thm:ci}]
	Without loss of generality we can suppose that $I = \{1, \dots, \ell\}$ for some $\ell < d$. By Lemma \ref{lem:Normal_Form}, there exists a matrix $A$ in Hermite normal form such that the support of $X$ is the lattice $L = A.\Ocal_K^d$. Then there exist a standard Gaussian vector $Y$ such that $X = AY$. Let $B = (A_{ij})_{\ell+1 \leq i,j \leq d}$ be the lower-right $(d-\ell)\times (d-\ell)$ block of $A$, and $\{f_i: \ell+1 \leq i \leq d\}$ the linear forms defined by the rows of $B$. For a subset $J := \{ j_1, \dots , j_r \}$ of ${1, \dots, d} \setminus I$, by Lemma \ref{lem:indep}, we have that $X_{j_1}, \dots, X_{j_r}$ are mutually independent given $X_I$ if and only if $\{f_{j_i}: 1 \leq i \leq r\}$ are orthogonal. Let us define the matrix $C$ as the matrix obtained from $B \in \GL_{d-\ell}(K)$ by scaling its rows so that they have norm $1$. Then from Proposition \ref{prop:orthogonality} we deduce that $X_{j_1}, \dots, X_{j_r}$ are independent given $X_I$ if and only if the images (modulo $\pi$) of the rows of $C$ indexed by $J$ are linearly independent in $k^{d-\ell}$. So conditional independence given $X_I$ is encoded in an $k$-realizable matroid given by the images modulo $\pi$ of the rows of $C$. 
\end{proof}

\begin{example}
	Let $p$ be a prime number and Consider the Gaussian vector $X = (X_1, X_2, X_3,X_4)^\top$ in $\Q_p^4$ with distribution given by the following lattice 
	\[ 
	L = \begin{bmatrix}
	1 & 0 &0 &0 \\ 1& 1 &0&0\\ 1 &0 & p &0 \\ 1&p^{-1}&p^{-1}&p^2
	\end{bmatrix} \Z_p^4.
	\] 
	If $I = \{1\}$ then the matrices $B$ and $C$ from the proof of \cref{thm:ci} in this case are
	\[
	B = \begin{bmatrix}
	1 &0&0\\  0 & p &0 \\ p^{-1} & p^{-1}&p^2
	\end{bmatrix} \in \Z_p^{3 \times 3} 
	\quad \text{ and } \quad
	C = \begin{bmatrix}
	1 \ & 0 \ &0\\  0 \ & 1 \ & 0 \\ 1 \ & 1\  &0
	\end{bmatrix}  \in \F_p^{3 \times 3}.
	\]
	Then from the matrix $C$ we deduce the following conditional independence statements 
	\[
	X_2 \indep X_4 | X_1,  \quad X_3 \indep X_4 | X_1 \quad \text{ and }  \quad  X_3 \indep X_4 | X_1,
	\] 
	However the three variables $X_2, X_3$ and $X_4$ are not mutually independent conditioned $X_1$. Independence conditioned on $X_1$ is encoded in the matroid that arises from the rows of $C$.

\end{example}

In the case where $\Xcal$ spans a proper subspace $W_{\Xcal} := \mathrm{span}_K(\Xcal)$ of $K^d$ we can define a Haar measure $\lambda$ on $W_{\Xcal}$ and the likelihood function defined for every full rank lattice of $W_{\Xcal}$ as $\Lcal(\Xcal,L) = \prod\limits_{x \in \Xcal} \frac{1_{x_i \in L}}{\lambda(L)}$. The maximum likelihood estimate in this case is again $L_{\Xcal} := \mathrm{span}_{\Ocal_K}(\Xcal)$, and it is the minimal lattice with respect to inclusion amongst those that maximize the likelihood.

\subsection{Tropicalization  of $K$-Gaussians}\label{subsec:tropical} For any $x=(x_1, \dots, x_d) \in K^d$, we denote by $\val(x) = (\val(x_1), \dots, \val(x_d))$ the vector of component-wise valuations of $x$. We call $\val(x)$ the \emph{tropicalization} of $x \in \K^d$. If $\Prob$ is a  probability distribution on $K^d$, we  call its push-forward by the valuation map $\val$, the \emph{tropicalization of $\Prob$}. 

For $d = 1$, the tropicalizations of Gaussians distributions are shifted geometric distribution on $\Z$ \cite{tran2018tropical}. For $d \geq 2$, the tropicalization of non-degenerate Gaussians  yields interesting distributions on $\mathbb{Z}^d$. To see this, let $X$ be a non degenerate Gaussian vector supported on a lattice $L$ with Hermite normal form given by the matrix $A$. Let us denote by $\Prob$ the probability distribution of $X$. We recall that $\Prob$ is the uniform probability distribution supported on $L$. Let $V := \val(X) $ be the image of $X$ under valuation. 

The tail distribution function of $V$ is the function $Q_{L}$ defined as follows:
\[
Q_{L} (v) := \Prob( V \geq v ),  \quad \text{for all } v \in \Z^d.
\]
For a vector $v = (v_1, \dots, v_d) \in \Z^d$ let us denote by $\bm{\pi^{v}}$ the lattice in $K^d$ defined by $\bm{\pi^v} = \diag(\pi^{v_1}, \dots, \pi^{v_d}) \Ocal_K^d$. Notice that we can rewrite the  function $Q_L$ as follows:

\begin{equation}\label{eq:exponent_func}
Q_L(v) = \Prob(X \in \bm{\pi^v}) = \lambda( L \cap \bm{\pi^v}).
\end{equation}

\begin{lemma}
	There exists $\varphi_{L} : \Z^d \rightarrow \Z_{\geq  0}$ such that $\ Q_{L} (v) = q^{ - \varphi_{L}(v)}$ for all $v \in \Z^d$.
\end{lemma}

\begin{proof}
	Follows immediately from \cref{eq:exponent_func} and \cref{lem:measure_exp}.
\end{proof}

We now get to the proof of \cref{thm:d2} which explains how $\varphi_L$ depends on the lattice $L$.

\begin{proof}[Proof of Theorem \ref{thm:d2}]
	Let $X$ be a non-degenerate Gaussian random variable in $K^2$ whose lattice $L$ has Hermite normal form
	\[  
	L = \begin{bmatrix}   \pi^{a}    & 0\\   \pi^{c}x & \pi^{b} \end{bmatrix} \Ocal_K^2 
	\]
	where $a \in \Z$ and $c \leq b \in \Z$ are integers and $x$ is an element in $K$ with valuation zero. Let $V = (V_1,V_2)$ be the tropicalization of $X$. For v = $(v_1,v_2) \in \Z^2$, we have
	\begin{align*}
	\Prob(V \geq v) &= \Prob(V_1 \geq v_1) \ \Prob(V_2 \geq v_2 | V_1 \geq v_1 )\\
	&= q^{-\max(v_1 - a,0)} \Prob(V_2 \geq v_2 | V_1 \geq v_1).
	\end{align*}
	The conditional distribution of $ X | V_1 \geq v_1$ is the Gaussian distribution on $K^2$ given by the lattice $L_{v_1}$  where  
	\[ 
	L_{v_1} = \begin{bmatrix}  \pi^{a + \max(v_1 - a,0)}         & 0\\    \pi^{c + \max(v_1 - a,0)} x & \pi^{b} \end{bmatrix} \Ocal_K^2.
	\]
	So the probability $\Prob(V_2 \geq v_2 | V_1 \geq v_1 )$ is given by:
	
	$$ \Prob(V_2 \geq v_2 | V_1 \geq v_1 ) = q^{- \max(v_2 - \min(b , c + \max(v_1 - a,0)), 0 )}.$$
	Therefore,
	\begin{align*}
	\varphi_{L}(v) &= \max(v_2 - \min(b , c + \max(v_1 - a,0)), 0 ) + \max(v_1 - a,0) \\
	&= \max( 0 , v_1 - a ,  v_2 - c ,  v_1 +  v_2 - a - b)
	\end{align*}
	Then $\varphi_L$ is the restriction to lattice points of the tropical polynomial $P_{L}$ given by \eqref{eqn:p.lambda}, where $c_{00} = 0$, $c_{10} = a$, $m_{01} = c$ and $c_{11} = a + b$. Since $b \geq c$ (by defnition of Hermite normal form of $L$), it is easy to check that \eqref{eqn:supermodular} holds. Conversely, suppose that $\phi$ satisfies the hypothesis in \cref{thm:d2}. One can reparametrize the coefficients of $P$ to obtain the integers $a,b,c$ and thereby the lattice $L$. 
\end{proof}

\begin{figure}[H]
	\centering
	\begin{subfigure}[b]{0.35\textwidth}
		\begin{tikzpicture}[x  = {(2cm,0cm)}, y  = {(0cm,2cm)}, scale = 1]
		
		\definecolor{pointcolor}{rgb}{ 0,0,0 }
		\tikzstyle{linestyle} = [color=black, thick]
		
		\coordinate (p1) at (0,0);
		\coordinate (p2) at (0.7,0.7);

		\node at (p1) [inner sep=7pt, right, black] {\tiny{$(a,c)$}};
		\node at (0.7,0.8) [inner sep=7pt, right, black] {\tiny{$(a + b - c , b)$}};
		
		\filldraw[pointcolor] (p1) circle (3 pt);
		\filldraw[pointcolor] (p2) circle (3 pt);

		\coordinate (e1) at (-1,0);
		\coordinate (e2) at (0,-1);
		\coordinate (e3) at (0.7,1.7);
		\coordinate (e4) at (2.7,0.7);
		
		\draw[linestyle] (e4) -- (p2) -- (e3);
		\draw[linestyle] (e1) -- (p1) -- (e2);
		\draw[linestyle] (p1) -- (p2);
		
		\node at (1,-0.2) {${v_1 - c_{10}}$};
		\node at (-0.5,-0.5) {${0}$};
		\node at (-0.5,1) {${v_2 - c_{01}}$};
		\node at (1.6,1.3) {${v_1 + v_2 - c_{11}}$};

		\end{tikzpicture}
		\caption{Regions of linearity of $P_{L}$.}
		\label{fig:RegionOfLin}
	\end{subfigure}
	~ 
	\begin{subfigure}[b]{0.5\textwidth}
		\begin{tikzpicture}[x  = {(2cm,0cm)}, y  = {(0cm,2cm)}, scale = 1]
		
		\definecolor{pointcolor}{rgb}{ 0,0,0 }
		\tikzstyle{linestyle} = [color=black, thick]
		
		\coordinate (p0) at  (6.25,1.5);
		\coordinate (p1) at  (8.25,1.5);
		\coordinate (p2) at  (6.25,3.5);
		\coordinate (p12) at (8.25,3.5);

		\node at (6,1.25) {$ (2,0,0)$};
		\node at (8.5,1) {$ (1,1,0)$};
		\node at (6,3.75) {$(1,0,1)$};
		\node at (8.5,3.75) {$ (0,1,1)$};
		
		\node at (5,1) {$ $};

		\filldraw[pointcolor] (p1) circle (3 pt);
		\filldraw[pointcolor] (p2) circle (3 pt);
		\filldraw[pointcolor] (p12) circle (3 pt);
		\filldraw[pointcolor] (p0) circle (3 pt);
		
		\draw[linestyle] (p0) -- (p1) -- (p12) -- (p2) -- (p0);
		\draw[linestyle] (p1) -- (p2);

		\end{tikzpicture}
		\caption{Regular subdivision of the Newton polytope of $P_L$ induced by its coefficients.}
		\label{trig:trig_dim2}
		
	\end{subfigure}  
	\caption{Description of the polynomial $P_{L}$ in dimension $d = 2$. The function $\varphi_L$ is the restriction of this polynomial to $\mathbb{Z}^2$.}
	\label{fig:TropPlots1}
\end{figure}
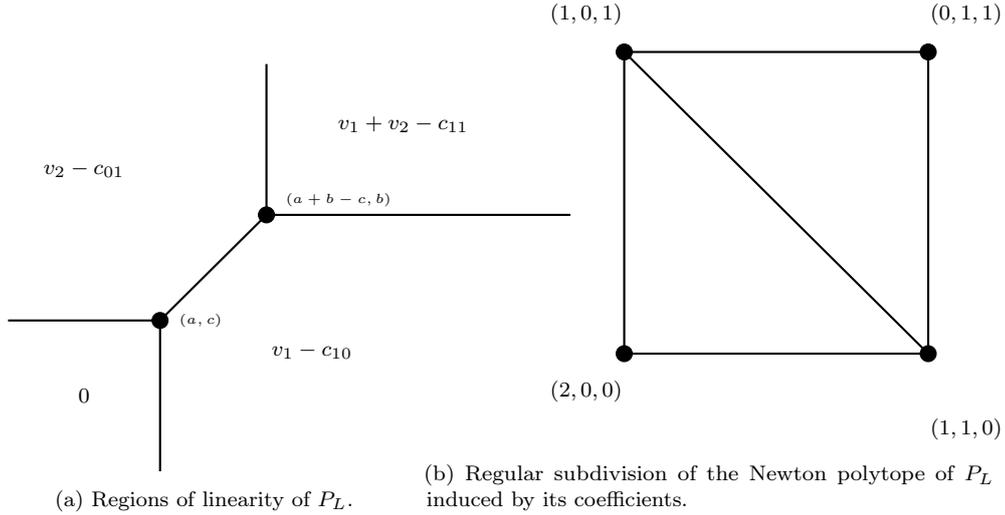

Note that, when $c = b$, the tropical variety of $P_{L}$ is the $1$-skeleton of the normal fan of a square. In that case the two entries of $X$ are independent and the probability distribution of $V$ is just a product of two shifted geometric distributions. When $c < b$, the rows of the lower triangular matrix are no longer orthogonal and thus the entries of $V$ are not independent by Lemma \ref{lem:indep}. This induces a unimodular triangulation on the square and the tropical variety of the polynomial $P_{L}$ takes the shape described in Figure \ref{fig:RegionOfLin}.

\section{Bruhat-Tits buildings}\label{sec:buildings}

For classical Gaussians, the cone of positive definite $ d \times d$ matrices is the parameter space for non-degenerate Gaussians on $\mathbb{R}^d$. For $K$-Gaussians, the role of covariance matrices is taken up by lattices in $K^d$, and the analogue of the positive semidefinite cone is the \emph{Bruhat-Tits building} for the group $\SL_d(K)$. In this section we shall briefly introduce some definitions of relevant material. For a reference on buildings, see \cite{abramenko2008buildings}. 

Let $L_1$ and $L_2$ be two full rank $K$-lattices in $K^d$. We say that $L_1$ and $L_2$ are equivalent if there exists a scalar $c \in K^{\times}$ such that $L_1 = c L_2$. This defines an equivalence relation on the set of full rank $K$-lattices in $K^d$. We denote by $[L]$ the equivalence class of a lattice $L$ i.e. its homotethy class. Two distinct equivalence classes $[L_1]$ and $[L_2]$ are said to be adjacent if there exists an integer $n \in \Z$ such that $\pi^{n+1} L_1 \subset L_2 \subset \pi^n L_1$. The rader can check that this relation is symmetric, but not transitive. The \emph{Bruhat-Tits building} $\Bcal_d(K)$ is the flag simplicial complex whose vertices are the equivalence classes of full rank lattices in $K^d$ and whose $1$-simplices are the adjacent equivalence classes. The simplices of $\Bcal_d(K)$ are the sets $\{ [L_1], \dots, [L_s]\}$ of distinct lattice classes that are pairwise adjacent. The following result relates adjacent cells of the Bruhat-Tits building with the matrix representations of the corresponding lattices.

\begin{proposition}\label{prop:adjacent}
	Let $S$ be the set elements $c.P$ such that $c \in K^{\times}$ and $P \in \GL_d(K) \cap  \Ocal_K^{d \times d} $ such that $\pi P^{-1} \in \Ocal_K^{d \times d}$ . Let $A,B \in \GL_d(K)$ be two invertible matrices and $L_A = A.\Ocal_K^d , L_B = B.\Ocal_K^d$ their associated lattices. Then the following hold
	
	\begin{enumerate}[label=(\roman*)]
		\item $[L_A] = [L_B]$ if and only if $A = cBU$ for some $c \in K^{\times}$ and $U \in  \GL_d(R)$.
		\item $[L_A]$ and $[L_B]$ are adjacent if and only if $A^{-1}B \in S$. 
	\end{enumerate}
	
\end{proposition}
\begin{proof}
	For $(i)$, suppose that $[L_A] = [L_B]$. Then, there exists $c \in K^\times$ such that $L_A = c L_B$. Then $A = cBU$ where $U$ is an orthogonal matrix in $\GL_d(R)$. Conversely if $A = cBU$ clearly the lattices $L_A $ and $L_B$ are in the same equivalence class. For $(ii)$, Suppose $[L_A]$ and $[L_B]$ are distinct adjacent classes. Then there exists $c \in K^\times$ such that $\pi L_A \subset c L_B \subset L_A$. Then there exist two invertible matrices $U$ and $V$ in $R^{d\times d }$ such that $\pi A = cBU$ and $cB = A V$. Then $VU = \pi I$ and $\pi V ^{-1} = U \in R^{d \times d} $. We have then $A^{-1} B = \frac{1}{c} V \in S$. Conversely we can easily obtain the inclusions $\pi L_A \subset c L_B \subset L_A$ when $A^{-1}B \in S$ holds.
\end{proof}

\begin{example} The building $B_2(\Q_p)$  is an infinite tree with degree $p+1$. \cref{fig:buildingb2} is a local depiction of $\Bcal_2(\Q_2)$.
	
	\begin{figure}[!ht]
		\centering
		\resizebox{0.6\textwidth}{!}{
			\begin{tikzpicture}
			\filldraw[] (0,0) circle [radius=0.1] node[above = 3.75mm] {$\begin{bmatrix}1 & 0 \\ 0 & 1\end{bmatrix}$};
			\draw[, thick] (0,0) -- (0, -4);
			\draw[, thick] (0,0) -- (3.464, 2);
			\draw[, thick] (0,0) -- (-3.464, 2);
			\filldraw[] (0,-4) circle [radius=0.1] node[above right = 1mm = 3.75mm] {$ \begin{bmatrix}1 & 0 \\ 1 & 2\end{bmatrix}$};
			\filldraw[]  (-1.732, -5) circle [radius = 0.1] node[below left=1.5mm] {$ \begin{bmatrix} 1 & 0 \\ 1 & 4\end{bmatrix}$};
			
			\draw[, thick] (0, -4) -- (-1.732, -5);
			\draw[, thick] (0, -4) -- (1.732, -5);
			\filldraw  (1.732, -5) circle [radius = 0.1] node[below right = 1 mm] {$\begin{bmatrix}1 & 0 \\ 3 & 4\end{bmatrix}$};
			
			\filldraw[] (3.464, 2)  circle [radius=0.1] node[below = 5mm] {$\begin{bmatrix}1 & 0 \\ 0 & 2\end{bmatrix}$};
			\filldraw[] (-3.464, 2) circle [radius=0.1] node[below = 4mm] {$\begin{bmatrix}2 & 0 \\ 0 & 1\end{bmatrix}$};
			\filldraw[] (3.464, 4) circle [radius = 0.1] node[right = 1.5mm]{$ \begin{bmatrix}1 & 0 \\ 2 &  4\end{bmatrix}$};
			\draw[, thick] (3.464, 2) -- (3.464, 4);
			\filldraw[] (5.2, 1) circle [radius = 0.1]      node[right = 1.5mm ] {$\begin{bmatrix}1 & 0 \\ 0 & 4\end{bmatrix}$};
			\draw[, thick] (3.464, 2) -- (5.2, 1);
			
			\filldraw[] (-3.464, 4) circle [radius = 0.1] node[right = 1.5mm]{$\begin{bmatrix}2 & 0 \\ 1 & 2\end{bmatrix}$};
			\draw[, thick] (-3.464, 2) -- (-3.464, 4);
			\filldraw (-5.2, 1) circle [radius = 0.1] node[left = 1.5mm]{$\begin{bmatrix}4 & 0 \\ 0 & 1\end{bmatrix}$};
			\draw[, thick] (-3.464, 2) -- (-5.2, 1);
			\end{tikzpicture}}
		\caption{}
		\label{fig:buildingb2}
	\end{figure}
	
\end{example}

Since the set $S$ in Proposition \ref{prop:adjacent} gives a certificate of adjacency for any pair of vertices, all vertices of $\Bcal_d$ have the same degree. We give a recipe to list the Hermite normal form of all lattice classes adjacent to the standard class $[\Ocal_K^d]$.

\begin{proposition}
	The $\Ocal_K$-modules $\Lambda$ satisfying $\pi \Ocal_K^d \subset \Lambda \subset \Ocal_K^d$ are exactly those with Hermite normal form $A = (a_{ij})$ satisfying the following conditions
	\begin{enumerate}[label =(\roman*)]
		\item $a_{ii} = \pi^{\epsilon_i}$ and $\epsilon_i = 0,1$ for all $1 \leq i \leq d$ and not all $\epsilon_i$ are $0$ or $1$,\\
		\item for $1 \leq j  < i \leq d$ we have the following:
		\begin{itemize}
			\item[(1)] 	$a_{ij} =0 $ if $a_{ii} =1 $ and $a_{i,j} \in \Tcal$ if $a_{ii} = \pi$,
			\item[(2)] 	$a_{ij} =0 $ if $a_{jj} =\pi $ and $a_{i,j} \in \Tcal$ if $a_{jj} = 1$.
		\end{itemize}

	\end{enumerate}
	In particular, the degree $deg(\Bcal_d(K))$ only depends on the dimension $d$ and the the module $q = \#k $ of the field $K$ and we have
	\[
	deg(\Bcal_d(K)) = \sum \limits_{l = 1}^{d-1}  \sum \limits_{1\leq k_1 < \dots < k_l \leq d} q^{k_1 + \dots + k_d - (1+ \dots +l)}.
	\] 
\end{proposition}

\begin{proof} We denote by $e_1, \dots, e_d$ the standard basis vectors of $K^d$. Let $C$ be a the homotethy class of full rank lattices which is adjacent to the standard lattices class $[\Ocal_K^d]$. Then, there exists a unique representative $L$ of $C$ such that $\pi \Ocal_K^d \subsetneq L \subsetneq \Ocal_K^d$. Let $A$ be the Hermite normal form of the lattice $L$. We consider $A_{ii} = \pi^{\epsilon_i}$ the diagonal coefficients of $A$ and $\epsilon = (\epsilon_1, \dots, \epsilon_d)$. Since $L \subset \Ocal_K^d$, we have $A \in \Ocal_K^{d \times d}$ so all entries of $A$ have non-negative valuations. The inclusion $p.\Ocal_K^d \subset L$ implies that all vectors of the form  $\pi e_i$  are $\Ocal_K$-linear combination of the columns of $A$. That means that $\epsilon_i \in\{0,1\}$. Notice that if all the $\epsilon_i$'s are equal to $0$, and since $A$ has non-negative valuation entries and is in Hermite normal form, we get $A = I_d$ is the identity matrix which contradicts $C \neq [\Ocal_K^d]$. Simillarly, if the $\epsilon_i$'s are all equal to $1$ then since $\pi \Ocal_K \subset L$ we deduce that $A = \pi I_d$ hence $L = p \Ocal_K^d$ and it is a contradiction. So, since $A$ has non-negative valuation entries and by definition of Hermite normal form, we deduce that $A$ satisfies the conditions $(i)$ and $(ii)-(1)$. It remains to show that the condition $(ii)-2$ is satisfied as well. Assume that there exists $1 \leq j < i \leq d$ such that $a_{jj} = \pi$ but $a_{ij} \neq 0$ and take $i$ to be minimal. Since $\pi e_j \in L$ and $L$ is generated over $\Ocal_K$ by the columns of $A$ we can write $\pi e_j$ as a $\Ocal_K$-linear combination of the columns of $A$. But since $A$ is lower triangular we can actually write $\pi e_j$ as the linear combination of the columns of $A$ indexed by $j, j+1, \dots, d$. So we can find an element $\alpha \in \Ocal_K$ such that $\alpha a_{i,i} = a_{ij} \neq 0$. So we deduce that $\val(a_{ii}) \leq \val(a_{ij})$ and this is a contradiction since $A$ is in Hermite normal form. So we deduce that $(ii)-(2)$ is also satisfied. Conversly, one can check without difficulty that if $A$ satifies conditions $(i)$ and $(ii)$ then $[L]$ is adjacent to $[\Ocal_K^d]$ in $\Bcal_d(K)$. Finally, the formula for the degree comes from simply counting the number of possibilities for the matrix $A$ under conditions $(i)$ and $(ii)$ which finishes the proof.
\end{proof}

\begin{example}
	Suppose $d=3$ and $K = \Q_p$ where $p$ is a prime number. Then the Hermite normal forms of the neighbours of the class $[\Ocal_K^3]$ is $\Bcal_d(\Q_p)$ are of the forms
	\begin{align*}
	\begin{bmatrix} 	1 & 0 & 0 \\	0 & 1 & 0 \\{\color{red}\ast} & {\color{red}\ast} & p \\ 	\end{bmatrix},	\quad   & 	\begin{bmatrix} 	1 & 0 & 0 \\	{\color{red}\ast} & p & 0 \\0 & 0 & 1 \\ 	\end{bmatrix}, \quad \begin{bmatrix} 	p & 0 & 0 \\	0 & 1 & 0 \\0 & 0 & 1 \\ 	\end{bmatrix}, \\
	\begin{bmatrix} 	1 & 0&0  \\	 {\color{red}\ast} & p & 0 \\{\color{red}\ast} & 0 & p \\ 	\end{bmatrix},	\quad   & 	\begin{bmatrix} 	p & 0 & 0 \\	0& p & 0 \\0 & 0 & 1 \\ 	\end{bmatrix}, \quad  \begin{bmatrix} 	p & 0 & 0 \\	0 & 1 & 0 \\0 & {\color{red}\ast} & p \\ 	\end{bmatrix}, \\
	\end{align*}
	where ${\color{red}\ast}$ represents any element in $\Tcal = \{0, \dots, p-1\}$. So, by this enumeration, the degree of $\Bcal_3(\Q_p)$ is
	\[
	\deg(\Bcal_3(\Q_p)) = 2 + 2p + 2p^2.
	\]	
\end{example}

The Bruhat-Tits building provides a pleasant geometric group theoretic framework for non-degenerate Gaussians on local fields. It is fruitful to tackle statistical problems with this in mind. 

Let $I = \{ 1, \dots, r\}$ and $J = \{ r+1, \dots ,r+s\}$ where $r,s \geq 0$ are integers with $r+s \leq d$. We denote by $\Mcal_{I,J}$ the model of non-degenerate Gaussian distributions on $\Q_p^d$ such that the variables indexed by $J$ are all independent given those indexed by $I$ i.e. the set of rank $d$ lattices $L$ in $\Q_p^{d}$ such that if $X$ is a Gaussian on with lattice $L$ we get
\begin{equation}\label{eq:indepStatement}
X_{r+1} \  \indep \ \dots \  \indep \ X_{r+s}  \ | \  X_1, \dots, X_r.
\end{equation}
Using \cref{lem:Normal_Form} and \cref{lem:indep} we can see that the points (or distributions) in $\Mcal_{I,J}$ are exactly those lattices $L$ whose Hermite normal form has the following shape
\[
\begin{bmatrix} 
\ast     & 0           & \dots       & \dots           &   \dots           & \dots     & \dots      & \dots                 & \dots        & 0 \\
\vdots & \ddots   & \ddots     &                    &                      &              &                &                           &                 &  \vdots\\
\ast     & \dots     & \ast         & \ddots        &                       &              &                &                          &                 &\vdots\\
\ast     & \dots     & \ast         & \ast            &   \ddots          &              &                &                          &                 &\vdots\\
\ast     & \dots     & \ast         & \nt{0}         &   \ast             & \ddots   &                 &                          &                 & \vdots\\
\vdots &              & \vdots      & \nt{\vdots} &   \nt{\ddots}      & \ddots   &  \ddots     &                          &                 &  \vdots\\
\vdots &              & \vdots      & \nt{0}         &   \nt{\dots}    & \nt{0}    &  \ast         &     \ddots          &                 & \vdots\\
\ast     & \dots    & \ast          & \ast            &   \dots           & \ast       &  \ast         &   \ddots            &   \ddots    & \vdots \\
\vdots &              & \vdots      & \vdots        &                      & \vdots   &  \vdots     &    \ddots              &   \ddots    & 0 \\
\ast     &  \dots   & \ast          & \ast            &        \dots      & \ast       &  \ast        &    \dots              &   \ast      & \ast \\
\end{bmatrix},
\]
where the red block is the below diagonal part of the block indexed by $J \times J$. The situation is similar to Gaussians on the real numbers where conditional independence implies that certain entries of the concentration matrix have to be zero.

The reader might be wary of the particular choice of $I$ and $J$, but by permuting the variables $X_1, \dots, X_d$ we may always assume that $I,J$ have the prescribed shapes. Also, since the conditional independence statement (\ref{eq:indepStatement}) does not change when scaling the lattice $L$ it is natural to want to describe the model $\Mcal_{I,J}$ in the Bruhat-Tits building.

\begin{problem}
	What does the model $\Mcal_{I,J}$ look like in the Bruhat-Tits building? What happens when we have multiple conditional independence statments? In the spirit of \cref{thm:mle}, can we easily fit conditional independence models to a collection of data points?
\end{problem}

\section{Summary and open questions}\label{sec:open}

This paper investigates statistical problems over local fields. We provided theorems on maximum likelihood estimation, conditional distribution, and distributions of tropicalized Gaussians. A major research question is to relate the building $\mathcal{B}_d(K)$ with statistical questions on Gaussians in $K^d$, such as conditional independence. Another direction is to generalize Theorem \ref{thm:d2} to $d \geq 3$. For a given lattice $L$, one can express it in canonical form and repeatedly condition on the values of $v_1$ and $v_2$ to compute an explicit expression for $\varphi_L$. We demonstrate in \cref{ex:d3}. Extensive computations in $K^3$ led us to Conjecture \ref{conj:d3} below.

\begin{example}\label{ex:d3}
	We consider the lattice 
	$  L = \begin{bmatrix}
	1    & 0         & 0    \\
	1    & \pi^2    & 0    \\
	1    & \pi      & \pi^2 
	\end{bmatrix} \Ocal_K^3  $. There exists a unique maximal (in the sens of inclusion) sublattice of $L$ representable by a diagonal matrix, we call this lattice the independence lattice of $L$, we denote it by $\ind(L)$ and in this case we have
	
	\[    \ind(L) = \begin{bmatrix}
	\pi^3    & 0     & 0    \\
	0     & \pi^3    & 0    \\
	0     & 0      & \pi^2 
	\end{bmatrix} \Ocal_K^3.
	\]
	We can compute the coefficients $c_I$ for $I \subset [3]$ and $|I| \leq 2$ using the proof of Theorem \ref{thm:d2} and all that is left is to compute the coefficient $c_{1,2,3}$. The computation of $ind(L)$ gives us the region of linearity of $P_L$ corresponding the monomial $c_{1,2,3} - v_1 - v_2 - v_3$ and in this case it is the orthant $\Ocal_{3,3,2} := \{ v \in \R^3 \colon v \geq (3,3,2)^T \}$. Using the coefficients we already computed we can then deduce $c_{1,2,3}$ and we find that
	\[
	P_{L} = \max(0, v_1, v_2, v_3, v_1 + v_2 - 2 , v_1 + v_3 - 1, v_2 + v_3 - 1, v_1 + v_2 + v_3 - 4).
	\]
	The support of $P_L$ is the unit cube $\{0,1\}^3$, with supermodular coefficients given by
	\[
	c_{\emptyset} = 0, c_{1} = 0, c_{2} = 0, c_{3} = 0, c_{12} = 2, c_{13} = 1, c_{23} = 1, c_{123} = 4.
	\]
	
	\begin{figure} [H]
		\centering
		
		\begin{subfigure}[b]{0.45\textwidth}
			\includegraphics[scale = 0.27]{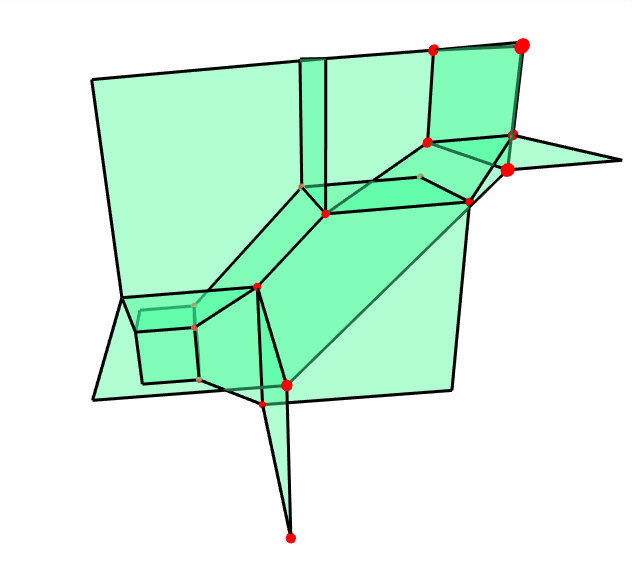}
			\caption{Tropical variety of $P_{L}$.}
		\end{subfigure}  
		~ 
		\begin{subfigure}[b]{0.45\textwidth}
			\includegraphics[scale = 0.5]{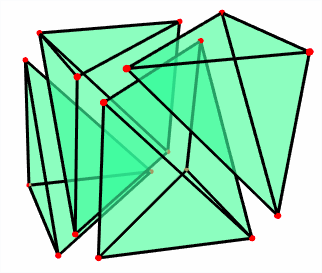}
			\caption{Regular subdivision of the Newton polytope of $P_L$ induced by its supermodular coefficients.}
		\end{subfigure}  
		
		\caption{ Tropical geometry of the Gaussian measure on $\lambda$}
		\label{Fig:tropical_Lattice}
	\end{figure}
	
\end{example}

\begin{conjecture}\label{conj:d3}
	Let $V$ be the tropicalization of a full-dimensional Gaussian in $K^d$ with lattice $L$. Define $\varphi_L: \mathbb{Z}^d \to \mathbb{R}$ via
	$$ \varphi_L(v) = - \log_q(\P(V \geq v)). $$
	Then a function $\phi: \mathbb{Z}^d \to \mathbb{R}$ equals to $\varphi_L$ for some $L$ if and only if $\phi$ is the restriction to lattice points of a tropical polynomial $P_L: \R^d \to \R$ supported on the $\{0,1\}^d$ cube with integer supermodular coefficients. 
	That is,
	\[
	P_{L}(v) = \max\limits_{I \subseteq [d]} \left(\sum_{i \in I} v_i - c_{I} \right),
	\] 
	where $(c_I)_{I \subset [d]}$ is a sequence of integers satisfying
	\[
		c_{\emptyset} = 0 \quad \text{ and } \quad  c_{I \cup J} + c_{I \cap J} \geq c_I + c_J, \quad \text{for all } I,J \subset \{1, \dots, d\}.
	\]
\end{conjecture}

\bibliography{article.bib}

\begin{thebibliography}{EVDD04}

\bibitem[AB08]{abramenko2008buildings}
Peter Abramenko and Kenneth~S Brown.
\newblock {\em Buildings: theory and applications}, volume 248.
\newblock Springer Science \& Business Media, 2008.

\bibitem[AGG12]{akian2012tropical}
Marianne Akian, Stephane Gaubert, and Alexander Guterman.
\newblock Tropical polyhedra are equivalent to mean payoff games.
\newblock {\em International Journal of Algebra and Computation},
  22(01):1250001, 2012.

\bibitem[AZ01]{MR1822891}
Sergio Albeverio and Xuelei Zhao.
\newblock A decomposition theorem for {L}\'evy processes on local fields.
\newblock {\em J. Theoret. Probab.}, 14(1):1--19, 2001.

\bibitem[BK13]{baldwin2013tropical}
Elizabeth Baldwin and Paul Klemperer.
\newblock Tropical geometry to analyse demand.
\newblock {\em Unpublished paper.[281]}, 2013.

\bibitem[Cas86]{cassels1986local}
John William~Scott Cassels.
\newblock {\em Local fields}, volume~3.
\newblock Cambridge University Press Cambridge, 1986.

\bibitem[CR66]{curtis1966representation}
Charles~W Curtis and Irving Reiner.
\newblock {\em Representation theory of finite groups and associative
  algebras}, volume 356.
\newblock American Mathematical Soc., 1966.

\bibitem[EL07]{evans2007expectation}
Steven~N Evans and Tye Lidman.
\newblock Expectation, conditional expectation and martingales in local fields.
\newblock {\em Electronic Journal of Probability}, 12(17):498--515, 2007.

\bibitem[Eva89]{evans1989local}
Steven~N Evans.
\newblock Local field gaussian measures.
\newblock In {\em Seminar on Stochastic Processes, 1988}, pages 121--160.
  Springer, 1989.

\bibitem[Eva93]{evans1993local}
Steven~N Evans.
\newblock Local field brownian motion.
\newblock {\em Journal of Theoretical Probability}, 6(4):817--850, 1993.

\bibitem[Eva01a]{evans2001local2}
Steven~N Evans.
\newblock Local field {U}-statistics.
\newblock {\em Contemporary Mathematics}, 287:75--82, 2001.

\bibitem[Eva01b]{evans2001local}
Steven~N Evans.
\newblock Local fields, gaussian measures, and brownian motions.
\newblock {\em Topics in probability and Lie groups: boundary theory},
  28:11--50, 2001.

\bibitem[Eva02]{evans2002elementary}
Steven~N Evans.
\newblock Elementary divisors and determinants of random matrices over a local
  field.
\newblock {\em Stochastic processes and their applications}, 102(1):89--102,
  2002.

\bibitem[Eva06]{MR2266718}
Steven~N. Evans.
\newblock The expected number of zeros of a random system of {$p$}-adic
  polynomials.
\newblock {\em Electron. Comm. Probab.}, 11:278--290, 2006.

\bibitem[EVDD04]{elsner2004max}
Ludwig Elsner and Pauline Van Den~Driessche.
\newblock Max-algebra and pairwise comparison matrices.
\newblock {\em Linear Algebra and its Applications}, 385:47--62, 2004.

\bibitem[Kac39]{kac1939characterization}
M~Kac.
\newblock On a characterization of the normal distribution.
\newblock {\em American Journal of Mathematics}, 61(3):726--728, 1939.

\bibitem[Khr13]{khrennikov2013p}
Andrei~Y Khrennikov.
\newblock {\em p-Adic valued distributions in mathematical physics}, volume
  309.
\newblock Springer Science \& Business Media, 2013.

\bibitem[LMY]{lin2018tropical}
Bo~Lin, Anthea Monod, and Ruriko Yoshida.
\newblock Tropical foundations for probability \& statistics on phylogenetic
  tree space.
\newblock {\em arXiv:1805.12400}.

\bibitem[Rob03]{ash2003numbertheory}
Ash Robert.
\newblock {\em A Course In Algebraic Number Theory}.
\newblock 2003.

\bibitem[Sch84]{schikhof1984ultrametric}
WH~Schikhof.
\newblock {\em Ultrametric Calculus (Cambridge Studies in Advanced Mathematics,
  4)}.
\newblock Cambridge University Press, Cambridge, 1984.

\bibitem[Sch07]{schikhof2007ultrametric}
Wilhelmus~Hendricus Schikhof.
\newblock {\em Ultrametric Calculus: an introduction to p-adic analysis},
  volume~4.
\newblock Cambridge University Press, 2007.

\bibitem[Tra13]{tran2013pairwise}
Ngoc~Mai Tran.
\newblock Pairwise ranking: Choice of method can produce arbitrarily different
  rank order.
\newblock {\em Linear Algebra and its Applications}, 438(3):1012--1024, 2013.

\bibitem[Tra18]{tran2018tropical}
Ngoc~Mai Tran.
\newblock Tropical gaussians: A brief survey.
\newblock {\em arXiv preprint arXiv:1808.10843}, 2018.

\bibitem[TY19]{tran2015product}
Ngoc~Mai Tran and Josephine Yu.
\newblock Product-mix auctions and tropical geometry.
\newblock {\em Math. Oper. Res.}, 44(4):1396--1411, 2019.

\bibitem[vR78]{van1978non}
Arnoud~CM van Rooij.
\newblock {\em Non-Archimedean functional analysis}.
\newblock Dekker New York, 1978.

\bibitem[VVZ94]{vladimirov1994p}
Vasilii~Sergeevich Vladimirov, Igor~Vasilievich Volovich, and Evgenii~Igorevich
  Zelenov.
\newblock {\em p-adic Analysis and Mathematical Physics}.
\newblock World Scientific, 1994.

\bibitem[Wei13]{weil2013basic}
Andr{\'e} Weil.
\newblock {\em Basic number theory.}, volume 144.
\newblock Springer Science \& Business Media, 2013.

\bibitem[YZZ19]{yoshida2017tropical}
Ruriko Yoshida, Leon Zhang, and Xu~Zhang.
\newblock Tropical principal component analysis and its application to
  phylogenetics.
\newblock {\em Bull. Math. Biol.}, 81(2):568--597, 2019.

\end{thebibliography}
\bibliographystyle{alpha}

\end{document}